\documentclass[psamsfonts,reqno]{amsart}
\usepackage{amssymb,eucal,graphics}

\newcommand{\two}{\mathbf{2}}

\newcommand{\bd}{\partial}
\newcommand{\pup}[1]{\textup{(}#1\textup{)}}
\newcommand{\sd}{\mathbin{\smallsetminus}}
\newcommand{\DD}{\mathbin{D}}
\newcommand{\ED}{\mathbin{\overline{D}}}
\newcommand{\tr}{\vartriangleleft}

\newcommand{\utr}{\trianglelefteq}

\newcommand{\dtr}{\mathbin{\vartriangleleft\kern-10pt {\lower
3pt\hbox{$\scriptscriptstyle\ne$}}\kern3pt}}

\newcommand{\jirr}{join-ir\-re\-duc\-i\-ble}

\newcommand{\jsd}{join-sem\-i\-dis\-trib\-u\-tive}
\newcommand{\jsdy}{join-sem\-i\-dis\-trib\-u\-tiv\-i\-ty}

\newcommand{\seq}[1]{\langle#1\rangle}

\newcommand{\set}[1]{\{#1\}}
\newcommand{\setm}[2]{\{#1\mid#2\}}

\newcommand{\at}{\mathsf{at}}
\newcommand{\es}{\varnothing}

\newcommand{\eps}{\varepsilon}
\newcommand{\SDJ}{\boldsymbol{SD}_{\vee}}

\DeclareMathOperator{\Max}{Max}
\DeclareMathOperator{\At}{At}

\newcommand{\SP}{\mathbf{S_p}}
\newcommand{\Co}{\mathbf{Co}}

\newcommand{\Subm}{\mathbf{Sub}_{\boldsymbol{\wedge}}}

\DeclareMathOperator{\J}{J}

\newcommand{\QQ}{\mathbb{Q}}

\newcommand{\KK}{\mathcal{K}}

\newcommand{\BIF}{\mathcal{BI}_{\mathrm{f}}}
\newcommand{\LB}{\mathcal{LB}}
\newcommand{\LBF}{\mathcal{LB}_{\mathrm{f}}}
\newcommand{\KF}{\mathcal{K}_{\mathrm{f}}}
\newcommand{\BI}{\mathcal{BI}}
\newcommand{\BS}{\mathbf{S}}
\newcommand{\BPO}{\mathbf{P}_{\boldsymbol{\omega}}}
\newcommand{\BPU}{\mathbf{P_u}}
\newcommand{\BQ}{\mathbf{Q}}

\numberwithin{equation}{section}

\theoremstyle{plain}

\newtheorem{lemma}{Lemma}[section]
\newtheorem{theorem}[lemma]{Theorem}
\newtheorem{proposition}[lemma]{Proposition}
\newtheorem{corollary}[lemma]{Corollary}

\newtheorem{claim}{Claim}

\newtheorem*{stat}{\name}
\newcommand{\name}{testing}

\theoremstyle{definition}
\newtheorem{definition}[lemma]{Definition}
\newtheorem{problem}{Problem}
\newtheorem*{notation}{Notation}

\theoremstyle{remark}
\newtheorem{remark}[lemma]{Remark}

\newcommand{\qedc}{{\qed}~{\rm Claim~{\theclaim}.}}
\newcommand{\qedsc}{{\qed}~{\rm Claim.}}

\newenvironment{cproof}
{\begin{proof}[Proof of Claim.]}
{\qedc\renewcommand{\qed}{}\end{proof}}

\begin{document}

\title[Finite biatomic lattices]
{Embedding finite lattices into\\
finite biatomic lattices}

\author[K.\,V.~Adaricheva]{Kira Adaricheva}
\curraddr{134 South Lombard Ave.\\
           Oak Park\\
           Illinois 60302\\
           USA}
\email{ki13ra@yahoo.com}

\author[F.~Wehrung]{Friedrich Wehrung}
\address{CNRS, UMR 6139\\
D\'epartement de Math\'ematiques\\
Universit\'e de Caen\\
14032 Caen Cedex\\
France}
\email{wehrung@math.unicaen.fr}
\urladdr{http://www.math.unicaen.fr/\~{}wehrung}
\date{\today}

\thanks{The second author was partially supported by the Fund of
Mobility of the Charles University (Prague), by
FRVS grant no. 2125, and by institutional grant
CEZ:J13/98:113200007.}


\begin{abstract}
For a class $\mathcal{C}$ of finite lattices, the question arises whether
any lattice in $\mathcal{C}$ can be embedded into some
atomistic, \emph{biatomic} lattice in $\mathcal{C}$. We provide
answers to the question above for $\mathcal{C}$ being, respectively,
\begin{itemize}
\item[---] The class of all finite lattices;

\item[---] The class of all finite \emph{lower bounded} lattices (solved by
the first author's earlier work).

\item[---] The class of all finite \emph{\jsd\ lattices} (this problem was,
until now, open).
\end{itemize}
We solve the latter problem by finding a quasi-identity valid in all
finite, atomistic, biatomic, \jsd\ lattices but not in all finite \jsd\
lattices.

\end{abstract}

\maketitle

\section{Introduction}
A lattice $L$ is \emph{biatomic}, if it is atomic (i.e., every
element of $L\setminus\set{0}$ lies above some atom of $L$) and whenever
$p$, $a$, and $b$ are elements of $L$ such that $p$ is an atom, $a$ and $b$
are nonzero, and $p\leq a\vee b$, there are atoms $x\leq a$ and $y\leq b$
such that $p\leq x\vee y$, see Definition~\ref{D:AtAtBiat}.

In our first result of this paper, Theorem~\ref{T:AllFinLatt}, we prove
that any finite lattice can be easily embedded atom-preservingly
into a finite biatomic one.

Biatomicity arises naturally in geometric lattices such as lattices of
subspaces of a vector space or, more generally, projective geometries. It
was also noticed by M.\,K.~Bennett \cite{Benn79} that in geometric
lattices, biatomicity is equivalent to modularity.

Biatomicity is probably even more common among \emph{convex geometries}.
The lattice theoretical facet of these at first finite and
purely combinatorial structures was studied in P.\,H.~Edelman \cite{Ed}
and P.\,H.~Edelman and R.~Jamison \cite{EdJa}.
In \cite{AGT}, these structures, now generalized to the infinite case,
were considered as an antithesis of geometric lattices, in terms of
the properties of the closure operators that define them.
We mention the lattices of convex subsets of a given affine space and
the lattice of subsemilattices of a given meet-semilattice as a few
examples of biatomic convex geometries.

Still, not all convex geometries are biatomic, thus to describe the biatomic
members within a given class of such structures would be of great interest.

Convex geometries are closely connected with the class of \jsd\ lattices.
A lattice $L$ is called \emph{\jsd}, if
  \[
  x \vee y = x \vee z \text{ implies that }
  x \vee y = x \vee (y\wedge z), \text{ for all } x,\,y,\,z \in L.
  \]
It is proved in \cite{AGT} that every finite \jsd\ lattice can be
embedded, in an atom-preserving way, into a finite, atomistic, \jsd\
lattice, or, equivalently, into a finite atomistic convex geometry.
This convex geometry is not generally biatomic.

The other construction in \cite{AGT} embeds any finite \jsd\ lattice
into the lattice $\SP(A)$ of algebraic subsets of some
(generally infinite) algebraic and dually algebraic lattice $A$.
This is also a convex geometry with the additional properties that it is
atomistic, biatomic, and \jsd.
This, together with \cite[Theorem~1.4]{AGT}, implies that \emph{every \jsd\
lattice $L$ can be embedded into an atomistic, biatomic, \jsd\ lattice $L'$},
see also the proof of \cite[Theorem~3.26]{AGT}.

It is asked in Problem~4 of \cite{AGT} whether $L'$ can be taken to
be finite whenever $L$ is finite. In the present paper, we solve this
problem in the negative, by showing a quasi-identity $\theta$ that is
satisfied by all finite, atomistic, biatomic,
\jsd\ lattices but not by all finite atomistic \jsd\ lattices, see
Theorem~\ref{T:BiatQuasid}. This result is inspired by a geometrical example
of finite convex geometry that in general produces non-biatomic \jsd\ lattices.

In contrast with this, we prove that every finite atomistic \jsd\ lattice $L$
can be $\seq{\vee,\wedge,0,1}$-embedded, in an atom-preserving way and with
the congruence extension property, into a finite, atomistic, \jsd\
lattice~$L'$ such that all biatomicity problems of~$L$ can be solved in $L'$,
see Theorem~\ref{T:PartSD+}.

We also study the case of finite
\emph{lower bounded} lattices, an important subclass of \jsd\ lattices.
The first author's earlier work \cite{Ad2} provides an embedding
of any finite lower bounded lattice into some finite biatomic convex geometry,
which implies Corollary~\ref{C:EmbSubP} (see also Theorem~\ref{T:PartSD+}).
Still, we do not know whether such an embedding can be done atom-preservingly.

This contributes to the list of open problems that concludes the paper.

\section{Biatomic lattices}\label{S:BasicBiat}

For a lattice $L$ with zero (i.e., least element), we denote by $\At L$
the set of atoms of $L$. The following definition recalls classical
notions, related to their counterparts in~\cite{Benn79}.

\begin{definition}\label{D:AtAtBiat}
A lattice $L$ with zero is
\begin{itemize}
\item \emph{atomic}, if every element of $L$ is above an atom of $L$;

\item \emph{atomistic}, if every element of $L$ is a join of atoms of
$L$;

\item \emph{biatomic}, if $L$ is atomic and for every atom $p$ of $L$
and all nonzero $a$, $b\in L$, if $p\leq a\vee b$, then there are atoms
$x\leq a$ and $y\leq b$ such that $p\leq x\vee y$.
\end{itemize}

\end{definition}

We observe that every finite lattice is atomic, and $L$ is atomistic
if{f} for all $a$, $b\in L$ such that $a\nleq b$, there exists
$p\in\At L$ such that $p\leq a$ and $p\nleq b$. The following lemma
is trivial.

\begin{lemma}\label{L:CharactBiat}
Let $L$ be an atomic lattice. Then $L$ is biatomic if{f} for every
atom $p$ of $L$ and all $a$, $b\in L\setminus\{0\}$ such that
$p\nleq a$, $p\nleq b$, and $p\leq a\vee b$, there exists an atom
$q\leq a$ of $L$ such that $p\leq q\vee b$.
\end{lemma}

For a lattice $L$ with zero and $x\in L$, let $\at(x)$ be the
statement that $x$ is an atom of $L$.
We can prove right away the following easy embedding result.

\begin{theorem}\label{T:AllFinLatt}
Let $L$ be a finite lattice. Then $L$ has a
$\seq{\vee,\wedge,0,1,\at}$-embedding into some finite,
atomistic, biatomic lattice $M$.
\end{theorem}

\begin{proof}
Put $A=L\setminus(\{0\}\cup\At L)$. For each $a\in A$, let $p_a$ and
$q_a$ be new distinct elements, and put
  \[
  M=L\cup\setm{p_a}{a\in A}\cup\setm{q_a}{a\in A}.
  \]
We define a partial ordering on $M$, extending the partial ordering of
$L$, by making all the elements of
$\setm{p_a}{a\in A}\cup\setm{q_a}{a\in A}$ mutually incomparable, and
by saying that
  \begin{gather*}
  x\leq p_a\text{ if{f} }x\leq q_a\text{ if{f} }x=0,\\
  p_a\leq x\text{ if{f} }q_a\leq x\text{ if{f} }a\leq x,
  \end{gather*}
for all $a\in A$ and $x\in L$. Then it is straightforward to verify that
$\leq$ is a lattice ordering on~$M$ and that the inclusion map from $L$
into $M$ is a $\seq{\vee,\wedge,0,1,\at}$-embedding. Since
$a=p_a\vee q_a$ for all $a\in A$, the lattice $M$ is atomistic.

To prove that $M$ is biatomic, it is convenient to use
Lemma~\ref{L:CharactBiat}. So let $p\in\At M$, let $a$, $b\in M$ such
that $p\nleq a$, $p\nleq b$, and $p\leq a\vee b$ (in particular, $a$
and $b$ are incomparable, thus they are nonzero), we find an atom
$q$ of $M$ such that $q\leq a$ and $p\leq q\vee b$. If $a$ is an atom of
$M$ then $q=a$ works, so suppose from now on that $a\in A$. If
$b\in L$, then $p\leq a\vee b=p_a\vee b$, thus $q=p_a$ is as required.
If $b\notin L$, say, $b=p_c$ for some $c\in A$ such that
$c\nleq a$, then $p\leq a\vee p_c=a\vee c=p_a\vee p_c$, so $q=p_a$ is
as required again.
\end{proof}

\section{Equivalence of definitions of embedding into finite biatomic
\jsd\ lattices}
\label{S:EquivDef}

We say that a partially ordered set is \emph{n{\oe}therian}, if it has no
infinite strictly increasing chain. Equivalently, every nonempty subset
has a maximal element. Of course, every finite partially ordered set is
n{\oe}therian.

In this paragraph we will prove that if an atomistic lattice can be embedded
into some n{\oe}therian, \jsd, biatomic lattice, then it also can be
embedded into an \emph{atomistic} such lattice. We will work toward
the proof of this statement in Corollary \ref{C:EquivBAE2}.

An immediate application of n{\oe}therianity gives the following
result.

\begin{lemma}\label{L:FinRep}
Let $L$ be a n{\oe}therian lattice and
let $G$ be a subset of $L$. Then every join of elements of $G$ is a finite
join of elements of $G$.
\end{lemma}

We recall the following elementary property of \jsd\ lattices that
follows immediately from \cite[Theorem~1.21]{FJN}.

\begin{lemma}\label{L:XcapY}
Let $K$ be a \jsd\ lattice with zero, 
let $a\in K$ and let $X$ and $Y$ be
finite sets of atoms of $K$. If $a\vee\bigvee X=a\vee\bigvee Y$, then
$a\vee\bigvee X=a\vee\bigvee(X\cap Y)$.
\end{lemma}

For a subset $S$ of a join-semilattice $L$, we denote the
set of all joins of nonempty finite subsets of $S$ by $S^\vee$.

\begin{proposition}\label{P:GenAt}
Let $L$ be a n{\oe}therian, biatomic, \pup{\jsd} lattice and let
$a\in L$. We put $S=[0,a]\cap\At L$. Then $T=(\set{0}\cup S)^\vee$ is a
n{\oe}therian, biatomic, \pup{\jsd}, atomistic lattice.
\end{proposition}

\begin{proof}
By definition, $T$ is a $\seq{\vee,0}$-subsemilattice of
$L$. Since $L$ is n{\oe}therian, 
it follows from Lemma \ref{L:FinRep} that $T$ is a lattice in its own right,
and it is n{\oe}therian. It is obviously atomistic, with $\At T=S$.

Let $p\in S$ and let $x$, $y\in T\setminus\set{0}$ such that
$p\leq x\vee y$. Since $L$ is biatomic, there are atoms $u\leq x$ and
$v\leq y$ of $L$ such that $p\leq u\vee v$. {}From $x$, $y\leq a$, it
follows that $u$, $v\in S$. This shows that $T$ is biatomic.

Now we prove that $T$ is \jsd\ provided $L$ is. Let $b$, $x$, $y\in T$ such
that $b\vee x=b\vee y$. 
By the definition of $T$, there are finite subsets $X$
and $Y$ of $S$ such that $x=\bigvee X$ and $y=\bigvee Y$. It follows from
Lemma~\ref{L:XcapY} that $b\vee x=b\vee\bigvee(X\cap Y)$. 
Since $\bigvee(X\cap Y) \leq x\wedge_Ty$, we get
$b\vee x=b\vee(x\wedge_Ty)$, thus completing the proof that $T$ is \jsd.
\end{proof}

For subsets $P$ and $Q$ in a lattice $L$, we say that $P$ \emph{separates
the elements of $Q$}, if for all $x$, $y\in Q$ with
$x\nleq y$, there exists $p\in P$ such that $p\leq x$ and
$p\nleq y$. Observe, in particular, that $L$ is atomistic if{f}
$\At L$ separates the elements of $L$.

\begin{corollary}\label{C:EquivBAE}
Let $L$ be a lattice. Then the following are equivalent:
\begin{enumerate}
\item There is a lattice embedding from
$L$ into some finite \pup{resp., n{\oe}therian}, \jsd, biatomic lattice
$M$ such that $\At M$ separates the elements of $L$.

\item There is a lattice embedding from $L$ into some finite
\pup{resp., n{\oe}therian}, \jsd, biatomic, atomistic lattice.
\end{enumerate}
The results above also hold for ``$0$-lattice embedding'' instead
of ``lattice embedding''.
\end{corollary}

\begin{proof}
We prove the nontrivial direction (i)$\Rightarrow$(ii), for
``n{\oe}therian''---the proof for ``finite'' is similar. Let $M$ be a
n{\oe}therian, biatomic, \jsd\ lattice containing $L$ as a sublattice
such that $\At M$ separates the elements of $L$. We put
  \[
  S=\setm{p\in\At M}{p\leq 1_L},
  \]
and we observe that $S$ separates the elements of $L$. Now
we put $T=(\set{0_M}\cup S)^\vee$. By Proposition~\ref{P:GenAt}, $T$ is a
n{\oe}therian, atomistic, biatomic, \jsd\ lattice, with $\At T=S$ and
$0_T=0_M$.

Now we shall investigate further the interaction of the lattices
$T$, $M$, $L$.

Define a map $f\colon L\to T$ by the rule
  \[
  f(x)=\bigvee\nolimits_T\setm{p\in S}{p\leq x},\quad\text{for all }x\in L.
  \]
In view of Lemma \ref{L:FinRep}, $f$ is well-defined.
We observe that $f(x)\leq x$, for all $x\in L$.

\begin{claim}\label{Cl:fmemb}
The map $f$ is a meet-embedding from $L$ into $T$. If $0_M=0_L$, then
$f(0_L)=0_T$.
\end{claim}

\begin{cproof}
It is clear that $f$ is order-preserving and that if $0_M=0_L$, then
$f(0_L)=0_T$. Let $x$, $y\in L$ such that $x\nleq y$. Since $S$
separates the elements of $L$, there exists $p\in S$ such that
$p\leq x$ and $p\nleq y$. Thus $p\leq f(x)$ (by the definition of $f$) and
$p\nleq f(y)$ (because $f(y)\leq y$), whence $f(x)\nleq f(y)$. So $f$
is an order-embedding.

Now let $x$, $y\in L$. We prove that
$f(x)\wedge_Tf(y)\leq f(x\wedge_Ly)$. Let $p\in S$ such that
$p\leq f(x)\wedge_Tf(y)$. Since $f(x)\leq x$ and $f(y)\leq y$, this
implies that $p\leq x$ and $p\leq y$. Thus, since $L$ is a sublattice
of $M$, the inequality $p\leq x\wedge_Ly$ holds, whence
$p\leq f(x\wedge_Ly)$. Therefore, since $f$ is order-preserving, $f$ is a
meet-homomorphism.
\end{cproof}

\begin{claim}\label{Cl:joinhom}
The map $f$ is a join-homomorphism from $L$ to $T$.
\end{claim}

\begin{cproof}
Let $x$, $y\in L$ and let $p\in S$ such that
$p\leq f(x\vee_Ly)$. Suppose that $p\nleq f(x)\vee_Tf(y)$. In
particular, $x$ and $y$ are nonzero in $L$, thus in $M$. By the definition
of $f$, this means that $p\leq x\vee_Ly=x\vee_My$. Hence, since $M$ is
biatomic and $x$ and $y$ are nonzero in $M$, there are atoms $u\leq x$ and
$v\leq y$ of $M$ such that $p\leq u\vee_Mv$. Observe that $u$, $v\in S$.
Moreover, $u\leq f(x)$ and $v\leq f(y)$, whence
$p\leq u\vee_Mv=u\vee_Tv\leq f(x)\vee_Tf(y)$, a contradiction.

Therefore, we have proved that $f(x\vee_Ly)\leq f(x)\vee_Tf(y)$. Since
$f$ is order-preserving, the converse inequality holds, which
concludes the proof of the claim.
\end{cproof}

The proof of Corollary~\ref{C:EquivBAE} is completed.
\end{proof}

\begin{corollary}\label{C:EquivBAE2}
Let $L$ be an atomistic lattice. Then the following are equivalent:
\begin{enumerate}
\item $L$ has a $\seq{\vee,\wedge,0}$-embedding into some finite
\pup{resp., n{\oe}therian}, \jsd, biatomic lattice.

\item $L$ has a $\seq{\vee,\wedge,0}$-embedding into some finite
\pup{resp., n{\oe}therian}, \jsd, biatomic, atomistic lattice.
\end{enumerate}
\end{corollary}

\begin{proof}
We prove the nontrivial direction (i)$\Rightarrow$(ii), for
``n{\oe}therian'' (the proof for ``finite'' is similar). Let $M$ be a
n{\oe}therian, \jsd, biatomic lattice
such that $L$ is a $0$-sublattice of
$M$. By Corollary~\ref{C:EquivBAE}, it is sufficient to prove that the
atoms of $M$ separate the elements of~$L$. So let $x$, $y\in L$ such that
$x\nleq y$. Since $L$ is atomistic, there exists $p\in\At L$ such that
$p\leq x$ and $p\nleq y$. Since $M$ is atomic, there exists an atom $q$
of $M$ below $p$. Then $q\leq x$. Furthermore,
$0\leq q\wedge y\leq p\wedge y=0$, and hence
$q\wedge y=0<q$, that is, $q\nleq y$. This proves our assertion.
\end{proof}

\section{Embedding finite lower bounded lattices}\label{S:LB}

For lattices $K$ and $L$, a lattice homomorphism $f\colon K\to L$ is
\emph{lower bounded}, if the preimage under $f$ of any principal dual ideal
of $L$ is either empty or has a least element. A lattice~$L$ is \emph{lower
bounded}, if every homomorphism from a finitely generated lattice to $L$ is
lower bounded. We refer the reader to \cite{AdGo,FJN} for more details.

For a finite meet-semilattice $P$, we denote by $\Subm(P)$ the lattice
of all subsemilattices of $P$ ($\es$ included). We state here the
following result from the first author's paper \cite{Ad2}.

\begin{theorem}\label{T:EmbSubP}
A finite lattice $L$ is lower bounded if{f} it can be embedded into
$\Subm(P)$ for some finite meet-semilattice $P$.
\end{theorem}

As $\Subm(P)$ is lower bounded, atomistic, and biatomic, this implies
immediately the following result.

\begin{corollary}\label{C:EmbSubP}
Any finite lower bounded lattice can be embedded
into some finite, atomistic, biatomic, lower bounded lattice.
\end{corollary}

For a finite, atomistic, lower bounded lattice $L$,
Theorem~\ref{T:EmbSubP} says that there exists an embedding from $L$
into $\Subm(P)$ for some finite meet-semilattice $P$. This embedding
can be chosen to preserve the zero, however, it may not preserve
atoms. The reason for this is that all lattices of the form $\Subm(P)$
have the property that for all atoms $x$ and $y$, there are at most
three atoms below $x\vee y$, while there are finite, atomistic, lower
bounded lattices that fail this property.

\section{One-atom extensions of finite atomistic lattices}\label{S:OneAt}

We start with the following definition.

\begin{definition}\label{D:ExtPair}
Let $L$ be a finite atomistic lattice. An \emph{extension pair} of $L$ is
a pair $(a;M)$, where the following are satisfied:
\begin{enumerate}
\item $a\in L\setminus(\set{0}\cup\At L)$;

\item $M$ is a meet-subsemilattice of $L$ that contains
$\set{0}\cup[a,1]$.
\end{enumerate}
For any $a\in L$, we put $L_a=L\setminus[a,1]$.
For an extension pair $(a;M)$, we put
  \[
  L(a;M)=(L_a\times\set{0})\cup(M\times\set{1}),
  \]
endowed with the componentwise ordering.
\end{definition}

By definition, a \emph{closure operator} of $L$ is a map $f\colon L\to L$
such that $f\circ f=f$, $f(x)\geq x$, and $x\leq y$ implies that
$f(x)\leq f(y)$, for all $x$, $y\in L$. Fix an extension pair $(a;M)$ of a
finite atomistic lattice $L$. Let $f$ be
the closure operator of $L$ associated with $M$, that is,
$f\colon L\to L$ is given by the rule
  \[
  f(x)=\text{least }y\in M\text{ such that }x\leq y,\qquad
  \text{for all }x\in L.
  \]
We observe that $L(a;M)$ is a meet-subsemilattice of $L\times\two$ (where
$\two=\set{0,1}$) that contains both $(0_L,0)$ and $(1_L,1)$ as elements.
Hence it is a lattice in its own right. For $(x,\eps)\in L\times\two$, we
denote by $\overline{(x,\eps)}$ the least element of $L(a;M)$ above
$(x,\eps)$. This element can easily be calculated, by the rule
  \begin{align*}
  \overline{(x,0)}&=\begin{cases}
  (x,0)&(\text{if }a\nleq x),\\
  (x,1)&(\text{if }a\leq x),
  \end{cases}\\
  \overline{(x,1)}&=(f(x),1).
  \end{align*}
We leave to the reader the straightforward proof of the following lemma.

\begin{lemma}\label{L:BasicLaM}
Let $L$ be a finite atomistic lattice and let $(a;M)$ be an extension pair
of $L$. Then the lattice $L(a;M)$ is finite atomistic, and the map
$j\colon L\to L(a;M)$ defined by $j(x)=\overline{(x,0)}$ for all $x\in L$
is a $\seq{\vee,\wedge,0,1,\at}$-embedding from $L$ into $L(a;M)$.
Furthermore, $\At(L(a;M))=(\At L\times\set{0})\cup\set{(0,1)}$.
\end{lemma}

Hence, $L(a;M)$ is an atomistic extension of $L$ by exactly one atom,
here $(0,1)$. In the sequel, the only properties of $L'=L(a;M)=L[p^*]$,
where $p^*=(0,1)$ is the new atom, that will be used are the ones listed
below:
  \begin{gather}
  L'=L\cup\setm{p^*\vee x}{x\in M}=
  L\cup\setm{p^*\vee x}{x\in L},\label{Eq:DescL'}\\
  p^*\leq x\Leftrightarrow a\leq x,\label{Eq:pleqx}\\
  x\leq p^*\vee y\Leftrightarrow x\leq f(y),\label{Eq:pleqxveey}
  \end{gather}
for all $x$, $y\in L$. Furthermore, $\At L'=\At L\cup\set{p^*}$. {}From now
on we shall use the more wieldy description of $L(a;M)$ given by
\eqref{Eq:DescL'}, \eqref{Eq:pleqx}, and \eqref{Eq:pleqxveey}.

\begin{remark}
It is not difficult to verify that conversely, every
$\seq{\vee,0,1}$-extension of~$L$ by exactly one atom below $1$ is, up to
isomorphism above $L$, of the form $L(a;M)$ for exactly one extension
pair $(a;M)$ of $L$. However, we shall not need this fact.
\end{remark}

Our next result describes when $L(a;M)$ is \jsd. For a subset $X$ of $L$,
we denote by $\Max X$ the set of all maximal elements of $X$.

\begin{lemma}\label{L:Lamjsd}
Let $L$ be a finite, atomistic, \jsd\ lattice and let $(a;M)$ be an
extension pair of $L$ with associated closure operator $f$. Then
$L(a;M)$ is \jsd\ if{f} the following conditions are satisfied:
\begin{enumerate}
\item $\Max L_a\subseteq M$;

\item $f(x\vee u)=f(x\vee v)$ implies that $u\leq f(x)$, for all $x\in L$
and all distinct atoms $u$ and $v$ of $L$.
\end{enumerate}
\end{lemma}

\begin{proof}
We put $L'=L(a;M)$. Suppose first that $L'$ is \jsd. Let $x\in\Max L_a$.
Suppose that $x\notin M$. Since $L$ is atomistic, there exists an atom
$u$ of $L$ such that $u\leq f(x)$ while $u\nleq x$. {}From $x<x\vee u$ and
the maximality of $x$ in $L_a=L\setminus[a,1]$, it follows that
$a\leq x\vee u$, whence $x\vee u\in M$. So, since $x\vee u\leq f(x)$,
we obtain that $x\vee u=f(x)$. Moreover, $p^*\leq a\leq x\vee u=f(x)$,
thus $x\vee u=f(x)\vee p^*=x\vee p^*$, and so, by the \jsdy\ of $L'$,
$u\leq x$, a contradiction. Therefore, $\Max L_a\subseteq M$.

Now let $x\in L$ and $u$, $v$ be distinct atoms of $L$ such that
$f(x\vee u)=f(x\vee v)$. It follows from \eqref{Eq:pleqxveey} 
that $x\vee u\vee p^*=x\vee v\vee p^*$,
whence, by the \jsdy\ of~$L'$, $u\leq x\vee p^*$, and therefore, 
again by \eqref{Eq:pleqxveey}, $u\leq f(x)$.

Conversely, suppose that both conditions (i) and (ii) are satisfied. To
prove the \jsdy\ of $L'$, if suffices to prove that $x\vee u=x\vee v>x$
cannot happen, for all $x\in L'$ and all distinct atoms $u$ and
$v$ of $L'$ (see \cite[Lemma~1.2]{AGT}). Since $L$ is \jsd, this holds if $x$,
$u$, $v\in L$.

Now suppose that $x\in L$ and $v=p^*$, so $x\vee p^*=x\vee u>x$. Hence, by
using \eqref{Eq:pleqxveey}, we obtain that
$x\vee u\leq f(x)\leq f(x)\vee p^*=x\vee p^*=x\vee u$, thus
$p^*\leq f(x)$, that is, by \eqref{Eq:pleqx}, $a\leq f(x)$. Moreover,
$x\vee u=f(x)>x$, so $x\in L_a$. Hence there exists $y\in\Max L_a$ such
that $x\leq y$. By assumption, $y\in M$, and consequently $f(x)\leq y$, a
contradiction since $a\leq f(x)$ and $a\nleq y$.

Since $x\vee u=x\vee v>x$, the last case to consider is $x=y\vee p^*$
for some $y\in L$ (see \eqref{Eq:DescL'}). It follows again from
\eqref{Eq:pleqxveey} that $f(y\vee u)=f(y\vee v)$, and consequently, by
assumption, $u\leq f(y)$, and so $x\vee u=x$, a contradiction.
\end{proof}

\section{Partially biatomic extensions}\label{S:BallRolls}

By definition, a ``biatomicity problem'' in a lattice $L$ is a formal
expression of the form $p\leq a\vee b$, where $p\in\At L$, $a$,
$b\in L\setminus\set{0}$, and the inequality $p\leq a\vee b$ holds while
$p\nleq a$, $p\nleq b$. A \emph{solution} of the problem above in $L$
consists of atoms $x\leq a$ and $y\leq b$ of~$L$ such that
$p\leq x\vee y$. Recall that a lattice embedding $f\colon K\hookrightarrow L$
has the \emph{congruence extension property}, if every congruence of $K$ is
the inverse image under $f$ of some congruence of $L$.

The present section will be mainly devoted to proving the
following results.

\begin{theorem}\label{T:PartSD+}
Every finite, atomistic, \jsd\ \pup{resp., lower bounded} lattice $L$
admits a $\seq{\vee,\wedge,0,1,\at}$-embedding with the congruence
extension property into some finite, atomistic, \jsd\ \pup{resp., lower
bounded} lattice $L'$ such that all biatomicity problems in $L$ can be
solved in $L'$.
\end{theorem}

\begin{remark}
It will turn out that the embedding from $L$ into $L'$ in
Theorem~\ref{T:PartSD+} preserves more than the congruences, it is in fact
an embedding for the transitive closure $\tr$ of the join-dependency
relation $\DD$. This is equivalent to $L'$ being a congruence-preserving
extension of $L$ in the finite, lower bounded case, but not in general.
\end{remark}

The core of the difficulty underlying Theorem~\ref{T:PartSD+} consists of
solving very special sorts of biatomicity problems. In Lemma~\ref{L:OneBB}
to Corollary~\ref{C:f,MCEP}, we let $L$ be a finite, atomistic, \jsd\
lattice, and $p$, $q$, $a\in L$ such that $p$ and $q$ are distinct atoms,
$a\in L\setminus(\set{0}\cup\At L)$, $p\leq a\vee q$, and
$p\nleq x\vee q$ for all $x<a$ in $L$.
Furthermore, we let $f\colon L\to L$ be the map defined by the rule
  \begin{equation}\label{Eq:TheClos}
  f(x)=\begin{cases}
  x&(\text{if }q\nleq p\vee x),\\
  p\vee x&(\text{if }q\leq p\vee x),
  \end{cases}
  \end{equation}
for all $x\in L$.

\begin{lemma}\label{L:OneBB}
The following assertions hold.
\begin{enumerate}
\item The map $f$ is a closure operator of $L$.

\item If we denote by $M$ the range of $f$, then $(a;M)$ is an extension
pair of $L$.

\item $L(a;M)$ is \jsd.

\item Denote by $p^*$ the unique atom of $L(a;M)\setminus L$. Then
$p<p^*\vee q$ and $p^*<a$.
\end{enumerate}
\end{lemma}

Hence, $L(a;M)$ is a \jsd\ extension of $L$ in which the biatomicity
problem $p\leq a\vee q$ has a solution.

\begin{proof}
The assertion (i) is straightforward. Furthermore, it is obvious that
$\set{0,1}$ is contained in $M$. Now let $x\in[a,1]$, we prove that $f(x)=x$.
This is obvious if $q\nleq p\vee x$, so suppose that $q\leq p\vee x$.
{}From $p\leq a\vee q$, $q\leq p\vee x$, and the \jsdy\ of $L$, it follows
that $p\leq x\vee a=x$, whence $f(x)=x\vee p=x$. This completes the proof
of (ii).

Now let $x\in\Max L_a$. We prove that $f(x)=x$. This is trivial if
$q\nleq p\vee x$, so suppose that $q\leq p\vee x$. If $p\nleq x$, then,
by the maximality assumption on $x$, $a\leq p\vee x$. Thus, from
$p\leq a\vee q$, $p\wedge a=0$, and the \jsdy\ of $L$, it follows that
$p\leq x\vee q$. Thus, since $q\leq p\vee x$ and by the \jsdy\ of~$L$,
we obtain that $p\leq x$, a contradiction. Therefore, $p\leq x$, so
$f(x)=p\vee x=x$. This proves that $\Max L_a\subseteq M$.

Let $x\in L$ and $u$, $v$ be distinct atoms of $L$ such that
$f(x\vee u)=f(x\vee v)$. We prove that $u\leq f(x)$.

If $q\nleq p\vee x\vee u$, then $q\nleq p\vee x\vee v$. Otherwise we would
have
  \[
  p\vee x\vee u=f(x\vee u)=f(x\vee v)=x\vee v,
  \]
thus $q\leq p\vee x\vee v=p\vee x\vee u$, a contradiction.
Thus $x\vee u=f(x\vee u)=f(x\vee v)=x\vee v$, whence, by the \jsdy\ of
$L$, $u\leq x\leq f(x)$.

Suppose now that $q\leq p\vee x\vee u$. By the previous paragraph,
$q\leq p\vee x\vee v$, and thus
$p\vee x\vee u=f(x\vee u)=f(x\vee v)=p\vee x\vee v$. Hence,
by the \jsdy\ of $L$, we have $u\leq p\vee x$, and so
$q\leq p\vee x\vee u=p\vee x$. Therefore, $u\leq p\vee x=f(x)$. By
Lemma~\ref{L:Lamjsd}, this completes the proof of assertion (iii).

The assertion (iv) follows immediately from $f(q)=p\vee q>p$.
\end{proof}

{}From Lemma~\ref{L:pDLu} to Corollary~\ref{C:f,MCEP}, we let $M$ and
$p^*$ be as in the statement and proof of Lemma~\ref{L:OneBB}. For a finite
lattice $K$, we let $\DD_K$ denote the relation of join-dependency on the set
of \jirr\ elements of $K$. Observe that for \emph{atoms} $x$ and $y$ of $K$,
the relation $\DD_K$ takes the following simple form:
  \[
  x\DD_K y\quad\text{if{f}}\quad x\neq y\text{ and }
  \exists u\in K\text{ such that }
  x\leq y\vee u\text{ and }x\nleq u.
  \]
Further, we denote by $\ED_K$ the binary relation on $\J(K)$ defined by
$x\ED_Ky$ if{f} either $x\DD_Ky$ or $x=y$. Then we let $\tr_K$ denote the
transitive closure of $\DD_K$ and $\utr_K$ denote the reflexive,
transitive closure of $\DD_K$.

Furthermore, since $L$ is finite, atomistic, and \jsd, it follows from
Lemma~\ref{L:XcapY} that every element $a$ of $L$
has a minimal decomposition, that is, a least (with respect to containment)
subset $X$ of $\At L$ such that $a=\bigvee X$. We denote this set of
atoms by $\bd^L(a)$ (``extreme boundary of $a$''), or $\bd(a)$ if $L$ is
understood. Note that $\bd(a)$ is also the unique irredundant decomposition
of~$a$. Observe that $\bd(a)$ consists
exactly of the elements which are join-prime in the interval $[0,a]$.
First it is convenient to prove the following lemma.

\begin{lemma}\label{L:pDLu}
For any $u\in\bd(a)$, the following relations hold:
\begin{enumerate}
\item $p\DD_Lu$;

\item $p^*\DD_{L[p^*]}u$.
\end{enumerate}
\end{lemma}

\begin{proof}
For any $u\in\bd(a)$, we put $a\sd u=\bigvee(\bd(a)\setminus\set{u})$.
{}From the fact that $u\in\bd(a)$, it follows that $a\sd u<a$, and so
$p\nleq(a\sd u)\vee q$ by the minimality assumption on $a$. However,
$p\leq a\vee q=(a\sd u)\vee u\vee q$ while $p\neq u$ (because $p\nleq a$),
whence $p\DD_Lu$.

Furthermore, $p^*\leq a=(a\sd u)\vee u$ while, since $a\sd u<a$,
we have  by \eqref{Eq:pleqx} that
$p^*\nleq a\sd u$, and consequently $p^*\DD_{L[p^*]}u$.
\end{proof}

\begin{lemma}\label{L:DrelSD}
For all $x$, $y\in\At L$, the following assertions hold:
\begin{enumerate}
\item $x\DD_{L[p^*]}y$ implies that $x\tr_Ly$;

\item $x\DD_{L[p^*]}p^*$ implies that $x\ED_Lp$;

\item $p^*\DD_{L[p^*]}x$ implies that there exists $u\in\bd(a)$ such that
$u\ED_Lx$;

\item $x\tr_{L[p^*]}y$ if{f} $x\tr_Ly$;

\item $p^*\tr_{L[p^*]}p^*$ if{f} there exists $u\in\bd(a)$ such that
$u\tr_Lp$.
\end{enumerate}
\end{lemma}

\begin{proof}
(i) By assumption, $x\neq y$ and there exists $u\in L[p^*]$ such that
$x\leq y\vee u$ and $x\nleq u$. Suppose that the relation $x\DD_Ly$ does
\emph{not} hold. So $u\notin L$, and therefore there exists $v\in M$ such
that $u=v\vee p^*$, hence, by \eqref{Eq:pleqxveey},
$x\leq f(y\vee v)$ and $x\nleq v$. Since the
relation $x\DD_Ly$ does not hold, we obtain that $x\nleq y\vee v$.
Hence $f(y\vee v)>y\vee v$, from where we obtain 
$f(y\vee v)=p\vee y\vee v$ (so $q\leq p\vee y\vee v$).
Then, since
$x\leq p\vee y\vee v$ and $x\nleq y\vee v$, we obtain that
  \begin{equation}\label{Eq:xEDLp}
  x\ED_Lp.
  \end{equation}
If $q\leq p\vee v$, then, since $v\in M$, the equalities $v=f(v)=p\vee v$
holds by the definition of $f$, whence $p\leq v$, and so $x\leq p\vee
y\vee v=y\vee v$, a contradiction. Hence $q\nleq p\vee v$, but $q\leq
p\vee y\vee v$, so we obtain the relation
  \begin{equation}\label{Eq:qEDLy}
  q\ED_Ly.
  \end{equation}
Finally, since $p\leq a\vee q$ and $p\nleq a$, we obtain that $p\DD_Lq$,
therefore, from \eqref{Eq:xEDLp} and \eqref{Eq:qEDLy}, it follows that
$x\tr_Ly$.

(ii) There exists $u\in L[p^*]$ such that $x\leq p^*\vee u$ and $x\nleq u$.
Thus $p^*\vee u\neq u$, so $u\in L$ and $x\leq f(u)$. {}From the relation
$x\nleq u$, it follows that $f(u)=p\vee u$, so $x\leq p\vee u$ while
$x\nleq u$, and so $x\ED_Lp$.

(iii) There exists $v\in L[p^*]$ such that $p^*\leq x\vee v$ and
$p^*\nleq v$. Thus $v\in L$, and $a\leq x\vee v$ while $a\nleq v$. {}From the
second relation, it follows that there exists $u\in\bd(a)$ such that
$u\nleq v$. However, $u\leq a\leq x\vee v$, and so $u\ED_Lx$.

(iv) {}From the fact that the natural embedding from $L$ into $L[p^*]$ is
atom-preserving, it follows that $x\tr_Ly$ implies that $x\tr_{L[p^*]}y$
for all $x$, $y\in\At L$. Conversely, for any $x$, $y\in\At L$, the
relation $x\tr_{L[p^*]}y$ means that there are a positive integer $n$ and
atoms $z_0=x$, $z_1, \dots, z_n=y$ of $L[p^*]$ such that
$z_i\DD_{L[p^*]}z_{i+1}$ for all $i<n$. We prove by induction on $n$ that
this implies that $x\tr_Ly$. For $n=1$, the conclusion follows from item
(i) above. Suppose that $n\geq2$. If $z_{n-1}\neq p^*$, then it follows
from the induction hypothesis that
$x\tr_Lz_{n-1}$, while, by item (i) above, $z_{n-1}\tr_Ly$, so $x\tr_Ly$.
Suppose now that $z_{n-1}=p^*$. Then $z_{n-2}\neq p^*$. Thus, by the
induction hypothesis, $x\utr_Lz_{n-2}$ (the equality may hold, \emph{e.g},
for $n=2$). Furthermore, it follows from items (ii) and (iii) above that
$z_{n-2}\ED_Lp$ and $u\ED_Ly$ for some $u\in\bd(a)$. But from
Lemma~\ref{L:pDLu}(i), it follows that $p\DD_Lu$, and so $z_{n-2}\tr_Ly$.
Therefore, $x\tr_Ly$.

(v) There exists $z\in\At L$ such that $p^*\tr_{L[p^*]}z\tr_{L[p^*]}p^*$.
{}From (ii), (iii), and (iv), it follows that $u\utr_Lz$,
for some $u\in\bd(a)$, and $z\utr_Lp$,
whence $u\utr_Lp$, but $u\neq p$ (because $p\nleq a$), and so $u\tr_Lp$.

Conversely, let $u\in\bd(a)$ such that $u\tr_Lp$. Thus we also have
$u\tr_{L[p^*]}p$. Since $p\leq p^*\vee q$
by Lemma \ref{L:OneBB} (iv)
and since $p$, $p^*$, and $q$ are distinct atoms, the relation
$p\DD_{L[p^*]}p^*$ holds. {}From Lemma~\ref{L:pDLu}(ii), it follows that
$p^*\DD_{L[p^*]}u$, so $p^*\DD_{L[p^*]}u\tr_{L[p^*]}p\DD_{L[p^*]}p^*$,
whence $p^*\tr_{L[p^*]}p^*$.
\end{proof}

\begin{corollary}\label{C:f,MCEP}\hfill
\begin{enumerate}
\item The canonical embedding from $L$ into $L[p^*]$ has the congruence
extension property; in fact, it is an embedding for the $\tr$ relation on
atoms.

\item If $L$ is lower bounded, then $L[p^*]$ is lower bounded.
\end{enumerate}
\end{corollary}

\begin{proof}
(i) By Theorem~2.30 and Lemma~2.36 in \cite{FJN}, it is sufficient to prove
that $x\utr_Ly$ if{f} $x\utr_{L[p^*]}y$, for all atoms $x$ and $y$ of~$L$,
which follows immediately from the stronger statement Lemma~\ref{L:DrelSD}(iv).

(ii) It is well-known that a finite lattice $K$ is lower bounded if{f} it has
no $D_K$-cycle, that is, the relation $\tr_K$ is \emph{irreflexive}, see
\cite[Corollary 2.39]{FJN}. Suppose that $L$ is lower bounded. It follows from
Lemma~\ref{L:DrelSD}(iv) that the relation $x\tr_{L[p^*]}x$ holds for no
$x\in\At L$. Suppose that $p^*\tr_{L[p^*]}p^*$. It follows from
Lemma~\ref{L:DrelSD}(v) that there exists $u\in\bd(a)$ such that
$u\tr_Lp$. By Lemma~\ref{L:pDLu}(i), $p\tr_Lu$, whence $L$ has a
$\DD_L$-cycle, a contradiction. Therefore, the relation $x\tr_{L[p^*]}x$
holds for no atom $x$ of $L[p^*]$.
\end{proof}

\begin{proof}[Proof of Theorem~\textup{\ref{T:PartSD+}}]
We present the proof for ``\jsd'', the proof for ``lower bounded'' is similar.
Since $L$ is finite, it suffices to prove that every biatomicity
problem $p\leq a\vee b$ in $L$ can be solved in some finite, atomistic,
\jsd\ $\seq{\vee,\wedge,0,1,\at}$-extension of $L$ in which $L$ has the
congruence extension property. We argue by induction
on $\ell_L(a)+\ell_L(b)$, where $\ell_L(x)$ denotes the
minimal size of a subset $X$ of $\At L$ such that $x=\bigvee X$, for all
$x\in L$. If $\ell_L(a)=\ell_L(b)=1$ then the biatomicity problem
$p\leq a\vee b$ is already solved in $L$, by $x=a$ and $y=b$. 
Now suppose, for example, that $b=c\vee q$, for some $c\in L\setminus\set{0}$
and some atom $q$ such that $\ell_L(c)<\ell_L(b)$.
Let $\overline{a}\leq a\vee c$ be minimal such
that $p\leq\overline{a}\vee q$.
By Lemma~\ref{L:OneBB}, there exists a finite \jsd\
$\seq{\vee,\wedge,0,1,\at}$-extension $L_1$ of~$L$, in which~$L$ has the
congruence extension property, such that there exists
an atom $p'\leq\overline{a}$ with $p\leq p'\vee q$. So $p'\leq a\vee c$
in~$L_1$ and
$\ell_{L_1}(a)+\ell_{L_1}(c)\leq\ell_L(a)+\ell_L(c)<\ell_L(a)+\ell_L(b)$.
Thus, arguing as above, we obtain a finite \jsd\
$\seq{\vee,\wedge,0,1,\at}$-extension $L_2$ of $L_1$, in which $L_1$ has the
congruence extension property, with atoms $x\leq a$
and $v\leq c$ such that $p'\leq x\vee v$. So
$p\leq p'\vee q\leq x\vee(v\vee q)$. Thus, again by Lemma~\ref{L:OneBB},
there exists a finite \jsd\
$\seq{\vee,\wedge,0,1,\at}$-extension $L_3$ of~$L_2$, in which $L_2$ has
the congruence extension property, with an atom
$y\leq v\vee q$ such that $p\leq x\vee y$. Observe that
$y\leq v\vee q\leq c\vee q=b$.
\end{proof}

\section{A quasi-identity for n{\oe}therian biatomic \jsd\ lattices}
\label{S:LjTr}

Let $\theta$ be the following quasi-identity in the language
$\seq{\vee,\wedge}$ of lattice theory:
\begin{multline*}
[u\leq a\vee b\vee v\ \&\
v\leq a\vee c\vee u\ \&\ (a\vee u)\wedge(b\vee c)\leq a\\
\ \&\ (a\vee b)\wedge(a\vee u)=(a\vee c)\wedge(a\vee v)=
(a\vee u)\wedge(a\vee v)=a]\\
\Longrightarrow u\leq a.
\end{multline*}
The present section will be mainly devoted to proving the
following result.

\begin{theorem}\label{T:BiatQuasid}
Every n{\oe}therian, atomistic, biatomic \jsd\ lattice with zero satisfies
$\theta$.
\end{theorem}

Let $M$ be a n{\oe}therian, atomistic, biatomic \jsd\ lattice with zero.
Observe that $M$ is a \emph{complete} lattice.
Let $a$, $b$, $c$, $u$, and $v$ be elements of $M$ satisfying the
premise of $\theta$, that is, the statement
  \begin{multline*}
  u\leq a\vee b\vee v\ \&\
  v\leq a\vee c\vee u\ \&\ (a\vee u)\wedge(b\vee c)\leq a\\
  \ \&\ (a\vee b)\wedge(a\vee u)=(a\vee c)\wedge(a\vee v)=
  (a\vee u)\wedge(a\vee v)=a.
  \end{multline*}
Suppose that $u\nleq a$, and put $S=\At M\setminus[0,a]$.
Since $M$ is atomistic, there exists $p\in S$ such that $p\leq u$.

\begin{lemma}\label{L:u0v0}
There are elements $u_0$, $v_0$ of $S^\vee$ such that the following
inequalities hold:
  \begin{equation}\label{Eq:u0v0}
  \begin{aligned}
  u_0\leq u&\text{ and }v_0\leq v;\\
  u_0\leq a\vee b\vee v_0&\text{ and }v_0\leq a\vee c\vee u_0.
  \end{aligned}
  \end{equation}
\end{lemma}

\begin{proof}
Suppose that $v\leq a$. Then $u\leq a\vee b$, thus
$u\leq(a\vee b)\wedge(a\vee u)=a$, a contradiction. Hence $v\nleq a$.

Put $x_0=p$; observe that $x_0\in S$.
So $x_0\leq u\leq a\vee b\vee v$, with $v$ nonzero (because
$v\nleq a$). Thus, since $M$ is
biatomic, there exists an atom $y_0$ of $M$ such that $y_0\leq v$ and
$x_0\leq a\vee b\vee y_0$. If $y_0\leq a$, then $x_0\leq a\vee b$, but
$x_0\leq u$, and so $x_0\leq u\wedge(a\vee b)\leq a$, a
contradiction. Hence $y_0\in S$.

Proceeding the same way with the inequality $y_0\leq a\vee c\vee u$ and
then inductively, we obtain elements $x_n$ and $y_n$, for $n<\omega$, of
$S$ such that $x_n\leq u$, $y_n\leq v$, $x_n\leq a\vee b\vee y_n$, and
$y_n\leq a\vee c\vee x_{n+1}$, for all $n<\omega$. Then
$u_0=\bigvee_{n<\omega}x_n$ and $v_0=\bigvee_{n<\omega}y_n$ (these are,
by Lemma~\ref{L:FinRep}, finite joins) are as required.
\end{proof}

Now, for $n<\omega$, suppose we have constructed $u_n$, $v_n\in S^\vee$
that satisfy the following inequalities:
  \begin{equation}\label{Eq:unvn}
  \begin{aligned}
  u_n&\leq a\vee b\vee v_n;\\
  v_n&\leq a\vee c\vee u_n;\\
  u_n&\leq a\vee u;\\
  v_n&\leq a\vee v.
  \end{aligned}
  \end{equation}
Since $u_n\leq b\vee(a\vee v_n)$ and $M$ is biatomic, for every
$x\in S\cap[0,u_n]$, there exists an atom $x^*\leq a\vee v_n$ of $M$
such that $x\leq b\vee x^*$. If $x^*\leq a$, then $x\leq a\vee b$.
However, $x\leq u_n\leq a\vee u$, and so
$x\leq(a\vee u)\wedge(a\vee b)=a$, a contradiction since $x\in S$. Hence,
$x^*\in S$, so that $v_{n+1}=\bigvee\setm{x^*}{x\in S\cap[0,u_n]}$ belongs
to $S^\vee\cap[0,a\vee v_n]$ and $u_n\leq b\vee v_{n+1}$. Proceeding in a
similar fashion with the inequality
$v_n\leq c\vee(a\vee u_n)$, we obtain elements $u_{n+1}$ and $v_{n+1}$ of
$S^\vee$ such that the following inequalities hold, see the
right half of Figure~1:
  \begin{equation}\label{Eq:uvn+1}
  \begin{aligned}
  u_{n+1}&\leq a\vee u_n;\\
  v_{n+1}&\leq a\vee v_n;\\
  u_n&\leq b\vee v_{n+1};\\
  v_n&\leq c\vee u_{n+1}.
  \end{aligned}
  \end{equation}
We verify that all the inequalities \eqref{Eq:unvn} are satisfied with $n$
replaced by $n+1$.
\begin{itemize}
\item $u_{n+1}\leq a\vee u_n\leq a\vee b\vee v_{n+1}$, and, similarly,
$v_{n+1}\leq a\vee c\vee u_{n+1}$.

\item $u_{n+1}\leq a\vee u_n\leq a\vee u$, and, similarly,
$v_{n+1}\leq a\vee v$.
\end{itemize}

Therefore, the values $u_0$ and $v_0$ obtained in Lemma~\ref{L:u0v0}
can be extended to sequences $(u_n)_{n<\omega}$ and $(v_n)_{n<\omega}$
of elements of $S^\vee$ that satisfy the inequalities listed in
\eqref{Eq:unvn} and \eqref{Eq:uvn+1} for all $n<\omega$.

A straightforward application of the last two inequalities in
\eqref{Eq:uvn+1} yields immediately the following lemma.

\begin{lemma}\label{L:bcuv}
The sequence $(b\vee c\vee u_{2n})_{n<\omega}$ is increasing.
\end{lemma}

Since $M$ is n{\oe}therian, there exists $n<\omega$ such that
$b\vee c\vee u_{2n}=b\vee c\vee u_{2n+2}$. Therefore, by using the last
two inequalities in \eqref{Eq:uvn+1}, we also obtain the following
equality:
  \begin{equation}\label{L:uv2n+1}
  b\vee c\vee u_{2n}=b\vee c\vee v_{2n+1}.
  \end{equation}

For any $n<\omega$, we let $U_n$ and $V_n$ be finite subsets of $S$ such
that $u_n=\bigvee U_n$ and $v_n=\bigvee V_n$. The existence of such sets
is ensured by Lemma~\ref{L:FinRep}.

\begin{lemma}\label{L:UVdisj}
$U_k\cap V_l=\es$, for all $k$, $l<\omega$.
\end{lemma}

\begin{proof}
Let $x\in U_k\cap V_l$. Then $x\leq u_k\leq a\vee u$ and
$x\leq v_l\leq a\vee v$, thus $x\leq(a\vee u)\wedge(a\vee v)=a$, which
contradicts the fact that $x\in S$.
\end{proof}

Now \eqref{L:uv2n+1} can be written as
$b\vee c\vee\bigvee U_{2n}=b\vee c\vee\bigvee V_{2n+1}$. But from
Lemma~\ref{L:XcapY} (applied to $K=M$) and Lemma~\ref{L:UVdisj}, it
follows that $U_{2n}\cap V_{2n+1}=\es$, and consequently
$u_{2n}=\bigvee U_{2n}\leq b\vee c$. However, $u_{2n}\leq a\vee u$, so
$u_{2n}\leq(a\vee u)\wedge(b\vee c)\leq a$, a contradiction since
$u_{2n}\in S^\vee$. This completes the proof of
Theorem~\ref{T:BiatQuasid}.

\begin{corollary}\label{C:LjTrnumber6}
There exists a finite, atomistic, \jsd\ lattice $L$ that cannot be
embedded into any finite \pup{or even n{\oe}therian} atomistic biatomic
\jsd\ lattice.
\end{corollary}

\begin{proof}
Put $L=\Co(\QQ^2,\set{a,b,c,u,v})$ where $a$, $b$, $c$, $u$, $v$ are as on
the left half of Figure~1, the lattice of all intersections with
$\set{a,b,c,u,v}$ of all convex subsets of $\QQ^2$, see \cite{AGT}. It is
well-known that all lattices of that form are \jsd.
  \begin{figure}[htb]
  \includegraphics{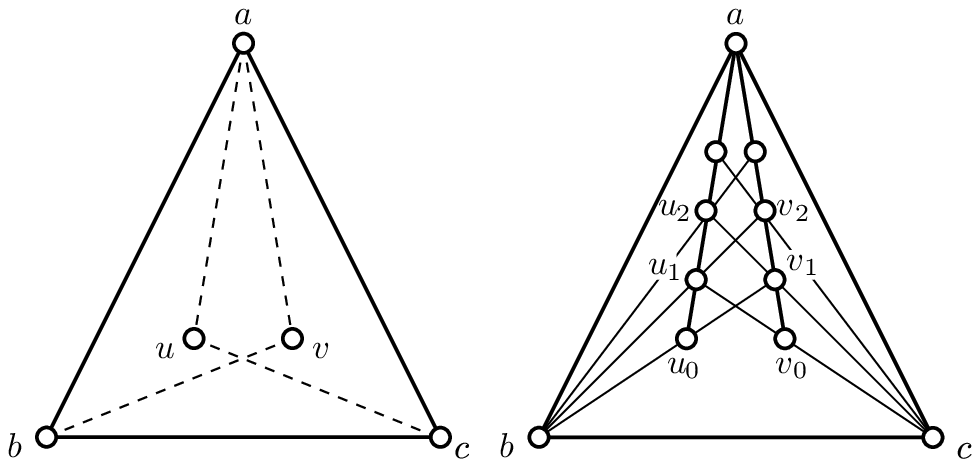}
  \caption{}
  \end{figure}
This configuration is obtained, for example, with
$a=\set{(0,3)}$, $b=\set{(-2,0)}$, $c=\set{(2,0)}$, $u=\set{(-1,1)}$, and
$v=\set{(1,1)}$. Then the premise of $\theta$ holds in $L$ for
these elements, although $u\nleq a$: hence $L$ does not satisfy
$\theta$.
\end{proof}

\begin{remark}
A closer look at the proof of Theorem~\ref{T:BiatQuasid} shows that
$L$ cannot be embedded into any finite, biatomic, \jsd\ lattice $M$ such
that there is an atom of $M$ below either $u$ or $v$ but not below $a$.
\end{remark}

\begin{remark}
The lattice $L$ can be embedded into an algebraic, atomistic, biatomic
convex geometry (see \cite{AGT}), namely, the lattice $\Co(\QQ^2)$ of all
convex subsets of $\QQ^2$.
In fact, $L$ can be embedded into a \emph{\jsd} atomistic biatomic
sublattice of $\Co(\QQ^2)$, namely, the lattice of convex polytopes of
$\QQ^2$, that is, finitely generated convex subsets of $\QQ^2$. Of
course, this lattice is neither n{\oe}therian, nor complete.

\end{remark}

\begin{remark}\label{Rk:DCC}
We recall that a partially ordered set is \emph{well-founded}, if every
nonempty subset has a minimal element. Then one can prove that the
quasi-identity $\theta$ is satisfied by every well-founded, biatomic,
\jsd\ lattice $L$, \emph{under the additional assumption that $u$ is an
atom of $L$}. We do not know whether the latter assumption can be
eliminated, see Problem~\ref{Pb:BiatQuasid}.
\end{remark}

\begin{remark}\label{RkLoCo}
It is proved in \cite{AGT} that every finite \jsd\ lattice has a
zero-preserving lattice embedding into a lattice of the form $\SP(A)$ (the
lattice of all algebraic subsets of $A$), for a lattice $A$ that is both
algebraic and dually algebraic. In particular, $\SP(A)$ is biatomic and
lower continuous. It follows that $\theta$ is not satisfied by all lower
continuous, atomistic, biatomic \jsd\ lattices.
\end{remark}

We observe the following immediate consequence of
Corollary~\ref{C:EquivBAE2} and Theorem~\ref{T:BiatQuasid}.

\begin{corollary}\label{C:LjTrNB0}
The lattice $L=\Co(\QQ^2,\set{a,b,c,u,v})$ of
Corollary~\textup{\ref{C:LjTrnumber6}} cannot be embedded into any
n{\oe}therian, biatomic, \jsd\ lattice.
\end{corollary}

\section{The quasivariety $\BQ(\BIF)$}\label{S:Quasiv}

\begin{notation}
We denote by $\BI$ the class of all atomistic, biatomic, \jsd\
lattices, and by $\LB$ the class of all lower bounded lattices.

We use standard notation for the basic operators defined on
$\seq{\vee,\wedge}$-structures, in particular, for a class $\KK$ of
$\seq{\vee,\wedge}$-structures, we define
\begin{itemize}
\item the class $\KF$ of all {\it finite} structures from $\KK$,
\item the class $\BS(\KK)$ of all structures that are embeddable into some
structure of $\KK$,

\item the class $\BPO(\KK)$ of all \emph{finite} direct products of structures
of $\KK$,

\item the class $\BPU(\KK)$ of all ultraproducts of structures of $\KK$,

\item the quasivariety $\BQ(\KK)$ generated by $\KK$.

\end{itemize}
\end{notation}
It was proved in \cite{AdGo} that $\LB\subset\BQ(\LBF)\subset\SDJ$,
both containments being proper. The class of biatomic \jsd\ lattices provides
a new element in this hierarchy. Our interest in this section will be focused
on the quasivariety generated by $\BIF$.

First we state that the finite members of this quasivariety are those
embeddable into lattices from $\BIF$.

\begin{proposition}\label{P:QBI}
The finite members of $\BQ(\BIF)$ are exactly the lattices that are
embeddable into some finite, atomistic, biatomic, \jsd\ lattice. In
formula, $\BQ(\BIF)_{\mathrm{f}}=\BS(\BIF)$.
\end{proposition}

\begin{proof}
It follows from results of the algebraic theory of quasivarieties
that the equality $\BQ(\KK)=\BS\BPU\BPO(\KK)$ holds for any class $\KK$, see,
for example, \cite[Corollary~2.3.4(3)]{Gorb}. We wish to prove that any
finite member $L$ of $\BQ(\BIF)$ is embeddable into some member
of~$\BIF$. Since the class $\BIF$ is closed under finite direct products,
that is, $\BPO(\BIF)\subseteq\BIF$, it follows from the
formula above that there exists a lattice embedding
$f\colon L\hookrightarrow L'$ where $L'\in\BPU(\BIF)$, that is, $L'$ is
an ultraproduct of members of $\BIF$.
Since $L$ is a finite system in a finite first-order language, a standard
argument about ultraproducts shows that $L$ can be embedded into some
system from~$\BIF$.
\end{proof}

Evidently, the proof above can be extended to any finite
first-order language, in particular the language $\seq{\vee,\wedge,0}$ if we
want to deal with lattices with zero, and so on.

\begin{proposition}\label{P:QBIF}
The following proper containments hold:

$\BQ(\LBF) \subset \BQ(\BIF) \subset \SDJ$.
\end{proposition}

\begin{proof}
It follows from Corollary~\ref{C:EmbSubP} that $\BQ(\LBF)\subseteq\BQ(\BIF)$.
Furthermore, the finite members of $\BQ(\LBF)$ are exactly the finite lower
bounded lattices while the lattice $\Co(\mathbf{4})$ of all 
order-convex subsets
of a four-element chain is finite, atomistic, biatomic, \jsd, and not lower
bounded, which shows that the containment above is proper.
The containment $\BQ(\BIF)\subseteq\SDJ$ holds by definition, and
Corollary~\ref{C:LjTrnumber6} provides an example of a finite \jsd\ lattice
which, by Proposition~\ref{P:QBI}, does not belong to $\BQ(\BIF)$.
\end{proof}

\section{Open problems}\label{S:Pbs}

According to Corollary~\ref{C:LjTrnumber6}, there exists a finite
atomistic \jsd\ lattice that cannot be embedded into any finite,
atomistic, biatomic, \jsd\ lattice. However, it is not hard to prove
that for all finite atomistic lattices $K$ and $L$ such that $K$
has a $\seq{\vee,0,\at}$-embedding into $L$, if $L$ is \jsd\
(resp. lower bounded), then so is $K$. Thus, in view of
Theorem~\ref{T:PartSD+}, the following question is natural.

\begin{problem}\label{Pb:vee0}
Let $L$ be a finite, atomistic, \jsd\ (resp., lower bounded)
lattice. Does $L$ have a $\seq{\vee,0,\at}$-embedding into some
finite, atomistic, biatomic, \jsd\ (resp., lower bounded) lattice?
\end{problem}

\begin{problem}\label{Pb:EmbDec}
For a finite \jsd\ lattice $L$, is it \emph{decidable} whether $L$ can be
embedded into some finite atomistic biatomic \jsd\ lattice?
\end{problem}

A variant of Problem~\ref{Pb:EmbDec} is the following.

\begin{problem}\label{Pb:FinAx}
Is the quasivariety $\BQ(\BIF)$ (see Section~\ref{S:Quasiv}) finitely 
based? That is, is the set of all quasi-identities satisfied by all finite
atomistic biatomic \jsd\ lattices equivalent to one of its finite subsets?
\end{problem}

By Proposition~\ref{P:QBI}, a positive answer to Problem~\ref{Pb:FinAx} would
imply a positive answer to Problem~\ref{Pb:EmbDec}. Nevertheless we
conjecture that Problem~\ref{Pb:FinAx} has a negative solution.

\begin{problem}\label{Pb:SDat}
Let $L$ be a finite lattice. If $L$ has a $\seq{\vee,\wedge,0}$-embedding
into some finite, atomistic, biatomic, \jsd\ lattice, then does $L$ have an
\emph{atom-preserving} such embedding?
\end{problem}

\begin{problem}\label{Pb:LBat}
Does every finite lower bounded lattice $L$ have an \emph{atom-preserving}
embedding into some finite, biatomic, lower bounded lattice?
\end{problem}

We have seen that every finite lower bounded lattice admits a
zero-preserving lattice embedding into some finite atomistic biatomic
lower bounded lattice, see Corollary~\ref{C:EmbSubP}.

Our final problem asks for extensions of Theorem~\ref{T:BiatQuasid}.

\begin{problem}\label{Pb:BiatQuasid}
Does any complete, \emph{upper continuous} (resp., \emph{well-founded}),
atomistic, biatomic, \jsd\ lattice satisfy the quasi-identity $\theta$
defined in Section~\ref{S:LjTr}?
\end{problem}

If we replace ``\jsd'' by ``convex geometry'' then the answer to the
corresponding problem is \emph{no}, as, for example, $\Co(\QQ^2)$ does
not satisfy $\theta$ while it is a complete, algebraic (thus upper
continuous), atomistic, biatomic convex geometry. However, $\Co(\QQ^2)$ is
not \jsd, see \cite[p.~234]{BB}. See also Remarks~\ref{Rk:DCC}
and~\ref{RkLoCo}.

\section*{Acknowledgments}
The work on this paper was initiated in February 2001 during both authors'
visit of Vanderbilt University, arranged by Prof. R.~McKenzie,
to whom we wish to express our gratitude.
The authors are grateful to the referees for careful reading of the paper.

\end{document}